\documentclass[11pt, english, reqno]{amsart}
\usepackage{fullpage}
\usepackage[utf8]{inputenc}
 \usepackage{amsmath,amssymb,stmaryrd,array,multirow,dsfont,mathrsfs}
\usepackage{amsthm}
\usepackage{enumerate}
\usepackage{xcolor}
\usepackage{color}
\usepackage{hyperref}
\usepackage{indentfirst}
\usepackage{mathabx}
\usepackage[utf8]{inputenc}
\usepackage[T1]{fontenc}
\usepackage{microtype}

\usepackage{graphics}
 \usepackage{amsmath,amssymb,stmaryrd,array,multirow,dsfont,mathrsfs}
\usepackage{amsthm}
\usepackage{xcolor}
\usepackage{color}
\usepackage{hyperref}

\graphicspath{{Images/}}

\xdefinecolor{darkgreen}{rgb}{0,0.4,0}

\newcommand{\tphi}{\widetilde{\phi}}
\newcommand{\tq}{\tilde{q}}

\newcommand{\beq}{\begin{equation}}
\newcommand{\eeq}{\end{equation}}
\newcommand{\beqs}{\begin{equation*}}
\newcommand{\eeqs}{\end{equation*}}
\newcommand{\ben}{\begin{eqnarray}}
\newcommand{\een}{\end{eqnarray}}
\newcommand{\beno}{\begin{eqnarray*}}
\newcommand{\eeno}{\end{eqnarray*}}
\renewcommand{\Re}{{\rm Re}\,}

\DeclareMathOperator{\Ran}{Ran}

\newcommand{\Rmnum}[1]{\uppercase\expandafter{\romannumeral #1} }
 \numberwithin{equation}{section}
\usepackage{mathtools}

\DeclarePairedDelimiterX{\inp}[2]{\langle}{\rangle}{#1, #2}

\newtheorem{thm}{Theorem}[section]

\newtheorem{prop}[thm]{Proposition}
\newtheorem{rmk}[thm]{Remark}

\def \d {\mathrm {d}}
\def\cA{{\mathcal A}}

\def\cC{{\mathcal C}}

\def\cL{{\mathcal L}}

\def\cO{{\mathcal O}}

\def\cR{{\mathcal R}}

\def\cU{{\mathcal U}}

\def\cW{{\mathcal W}}

\let\f=\frac
\def \p {\partial}
\def\mZ {\mathbb{Z}}
\def\mR {\mathbb{R}}
\def\mN {\mathbb{N}}

\def \pt {\partial_{t}}

\def \diag {\textnormal{diag}}

\def \uu {\underline{u}}

\def \uvd {\underline{v}_{\delta}}
\def \uv {\underline{v}}

\def \lxd {L_{\xi}^{\delta}}
\def \pxd {\Pi_{\xi}^{\delta}}

\title{Linear asymptotic stability of small-amplitude periodic waves of the generalized Korteweg--de Vries equations}

\author{Corentin Audiard}
\address{Sorbonne Université, CNRS, Université de Paris, Laboratoire Jacques-Louis Lions (LJLL), F-75005 Paris, France }
\email{{\tt corentin.audiard@upmc.fr}}
\thanks{Research of C.A. was partially supported by the French ANR Project NABUCO ANR-17-CE40-0025.}

\author{L.~Miguel Rodrigues}
\address{
Univ Rennes \& IUF, CNRS, IRMAR - UMR 6625, F-35000 Rennes, France}
\email{{\tt luis-miguel.rodrigues@univ-rennes1.fr}}
\thanks{Research of L.M.R. was partially supported by EPSRC grant no EP/R014604/1.}

\author{Changzhen Sun}
\address{Universit\'e de Toulouse, CNRS, IMT - UMR 5219, UPS, F-31062 Toulouse Cedex 9, France}
\email{changzhen.sun@math.univ-toulouse.fr}
\thanks{Work of  C.S. is supported by the ANR LabEx CIMI (grant ANR-11-LABX-0040) within the French State Programme ``Investissement d'Avenir''}

\date{\today}

\begin{document}

\begin{abstract}
We extend the detailed study of the linearized dynamics obtained for cnoidal waves of the Korteweg--de Vries equation in \cite{JFA-R} to small-amplitude periodic traveling waves of the generalized Korteweg--de Vries equations that are not subject to Benjamin--Feir instability. With the adapted notion of stability, this provides for such waves,  global-in-time bounded stability in any Sobolev space, and asymptotic stability of dispersive type. When doing so, we actually prove that such results also hold for waves of arbitrary amplitude satisfying a form of spectral stability designated here as dispersive spectral stability. 
\vspace{1em}


\vspace{1em}

\end{abstract}

\maketitle

\section{Introduction}

\subsection*{Outline}

In this note, we are concerned with the linear stability of 
periodic traveling waves  of the generalized  Korteweg-de Vries (gKdV) equation
\begin{align}\label{gkdv}
    \pt u+ \p_x(\p_x^2u+ f(u))=0\,.
\end{align}
where $f$ is a smooth function. 
This encompasses the power-type case when $f$ is a nonlinear monomial, that is, $f(u)=a\,u^p$, for some $p\in \mN$, $p\geq2$, and some $a\neq0$. The latter class includes the classical Korteweg--de Vries equation, obtained with $p=2$, and the modified Korteweg--de Vries equations, obtained with $p=3$, which are known to be integrable. 

A traveling wave of \eqref{gkdv}, of speed $c$, is a solution taking the form $(t,x)\mapsto \uu(x-ct)$ for some profile $\uu$. Such a wave is called periodic if its profile is a periodic function. A key point is that even if our background solutions are spatially periodic the dynamics we are interested in is set on the real line, that is, $x\in\mR$. In particular, in accordance with most of related concrete applications, we aim at understanding how a periodic traveling wave behaves under smooth localized perturbations. 

The study of the generalized Korteweg--de Vries equations and, even of the more specialized question of the stability of their traveling waves, have grown too considerably to be properly reported in the present note. We simply refer the reader to \cite{Zhidkov,Tao,Linares-Ponce,Angulo-Pava,Kapitula-Promislow} as possible entering gates to various aspects of the field. On the more specific question of periodic waves of such equations we refer to \cite{R,R_Roscoff} and the introductions of \cite{Nonlinearity-BMR2,NRS} for a general picture and a comprehensive list of references.

We only point out that, with the exception of the linear analysis in \cite{JFA-R} and some case-studies of integrable equations \cite{Mikikits-PhD,Mikikits-Teschl}, the stability analysis of periodic waves of dispersive equations is for the moment essentially restricted to spectral studies. Note moreover that though most of the analysis of \cite{JFA-R} relies on robust arguments, it does use at some key steps detailed spectral information obtained in \cite{Bottman-Deconinck} from the integrability through a Lax pair representation of the classical Korteweg--de Vries equation.

An inspection of \cite{JFA-R} shows that integrability is used there in two fundamentally different ways: on one hand to check a strong form of spectral stability that we designate as dispersive spectral stability, on the other hand to obtain a suitable form of high-frequency expansions of Bloch eigenfunctions. 

In the present note, motivated by this observation, we achieve two tasks. On one hand we replace the high-frequency analysis of \cite{JFA-R} with a robust one, hence proving that for any periodic wave of (gKdV) dispersive spectral stability implies linear asymptotic stability and comes with the detailed asymptotic behavior provided in \cite{JFA-R}. On the other hand we prove dispersive spectral stability for small-amplitude periodic waves that are not subject to Benjamin--Feir instability.

We regard the precise identification of the right notion of dispersive spectral stability and the high-frequency analysis as the main contributions of the present note. Indeed, though not written down elsewhere in the literature, a large part of the spectral stability analysis for small-amplitude waves could be considered, to some extent, as common knowledge of the community, or, in other words, as part of the folklore. In particular, the Benjamin--Feir instability criterion has been derived in a considerably more general framework in \cite{Nonlinearity-BMR2,SMF-Audiard-Rodrigues} (by combining the small-amplitude expansions of \cite{IUMJ-BMR} with the general spectral criterion of \cite{JNS-BNR,SMF-Audiard-Rodrigues}). Likewise, the fact that when the Benjamin-Feir instability is precluded spectral stability --- in the sense that the spectrum lies on the imaginary axis --- holds for small-amplitude waves has been proved for power-type nonlinearities in \cite{PhysD-Haragus-Kapitula} (with arguments close to those used in the related \cite{JDE-Gallay-Haragus}). In the stability direction our spectral analysis merely extends the standard argument to more general nonlinearities and, more importantly, strengthens the spectral conclusions so as to reach a stage when the linear analysis of \cite{JFA-R} completed with the present high-frequency analysis could be applied. 

\subsection*{Main results}  

To provide precise statements, let us fix a periodic wave of (gKdV) of profile $\uu$ and speed $c$. In a frame moving at speed $c$ this corresponds to a stationary solution and the dynamics linearized about it obeys $\p_t u=\cL u$ where $\cL:=-\p_x((\p_x^2+f'(\uu)-c)\,\cdot\,)$. When discussing spectral properties we regard $\cL$ as acting on $L^2(\mR)$ with domain $H^3(\mR)$ but we shall consider various extensions when turning to dynamical considerations.

It is convenient to use Floquet-Bloch theory to describe the spectrum of $\cL$. We recall some rudiments of it in the next section but for the moment it is sufficient to know that the spectrum of $\cL$ is given by a discrete set of curves locally parameterized by a real number called Floquet exponent $\xi$ and that multiplicity refers to multiplicity at a fixed Floquet exponent.

We say that the wave associated with $(\uu,c)$ is dispersively spectrally stable if the following list of properties holds. 
\begin{enumerate}
\item The spectrum of $\cL$ lies on the imaginary axis.
\item Each nonzero element of the spectrum is simple and the second-order derivative with respect to the Floquet exponent of the spectral curve passing through it is nonzero.
\item The spectral element $0$ corresponds to the Floquet exponent $0$ and has multiplicity three. Moreover the derivatives of the three curves passing through it are distinct and their third-order derivatives are nonzero. 
\end{enumerate}
The first property alone is referred to as spectral stability and, correspondingly, failure of it is known as spectral instability. Concerning the third one, we point out that it is a non trivial but classical fact of the periodic wave stability theory that the corresponding curves are always at least $\cC^1$ and that they are analytic when their first-order derivatives are distinct.

It follows from elementary considerations that $f'(\uu_0)-c_0\geq 0$ is necessary for the existence of a family of periodic waves with profiles $\uu$ converging to a constant value $\uu_0$ and speed $c$ converging to $c_0$, and that $f'(\uu_0)-c_0>0$ is sufficient to guarantee such an existence. With this in mind we may state our first main result.

\begin{thm}[Spectral stability]\label{thm1}
Let $(\uu_0,c_0)\in\mR^2$ satisfy $f'(\uu_0)-c_0>0$ and set
  \beq\label{cond-stability} 
\Delta_{BF}:=(f''(\uu_0))^2-3\,f'''(\uu_0)(f'(\uu_0)-c_0)\,.
  \eeq
\begin{enumerate}
\item If $\Delta_{BF}<0$, periodic waves associated with $(\uu,c)$ sufficiently close to $(\uu_0,c_0)$ are spectrally unstable.
\item If $\Delta_{BF}>0$, periodic waves associated with $(\uu,c)$ sufficiently close to $(\uu_0,c_0)$ are dispersively spectrally stable.
\end{enumerate}
\end{thm}
  
In order to analyze linear stability, we need to use the right way to think about proximity. With this goal in mind, a linear analogous to the notion of space-modulated stability introduced in \cite{Inventiones-JNRZ} was designed in \cite{JFA-R}. It consists in measuring the size of a function $u$ with
\beq\label{def-norm}
 N_{X}(u):=\inf_{\substack{(w,\psi)\\(u,\p_x\psi)\in X^2\\
 u=w+\psi \uu'}}(\|w\|_{X}+\|\p_x\psi\|_{X}),
\eeq
for some choice of a functional space $X$ endowed with a norm $\|\cdot\|_X$. We refer to \cite{R_Roscoff,JFA-R} for a detailed discussion of this choice. We only mention that if we were, instead, studying solitary waves or kinks, we would measure the size of $u$ with
\[
\inf_{\substack{(w,\psi_0)\in X\times \mR\\
 u=w+\psi_0\uu'}}\|w\|_{X}
\]
so as to obtain a linear analogous of the classical notion of orbital stability.

With this in hands, let us denote $(S_{\cL}(t))_{t\in\mR}$ the one-parameter group generated by $\cL$ and state our second main result. 

\begin{thm}[Linear stability]\label{thm2}
Assume that the wave associated with $(\uu,c)$ is dispersively spectrally stable.
\begin{enumerate}
\item For any $s\in \mR,$ there exists $C_s>0$ such that for any $t\in\mR$
\beqs 
  N_{H^s(\mR)}(S_{\cL}(t)(u_0))\leq C_s  N_{H^s(\mR)}(u_0).
\eeqs
\item There exits $C>0$ such that for any $t\in\mR$
  \beqs
  N_{L^{\infty}(\mR)}(S_{\cL}(t)(u_0))\leq C\,(1+|t|)^{-\f{1}{3}}\,N_{L^{1}(\mR)\cap H^1(\mR)}(u_0)\,.
  \eeqs
\end{enumerate}
\end{thm}

With the same assumption as in the foregoing theorem we may also provide a detailed description of the dynamics. Roughly speaking Theorem~\ref{thm2} proves that up to an error of order $|t|^{-1/3}$ the linearized dynamics is described by the arising of a space-time modulation of phase shift, $\psi(t,x)\,\uu'$ being thought as a linear part of $\uu(x+\psi(t,x))$. By adding a space-time modulation of the other parameters describing the wave one captures the dynamics up to a $|t|^{-1/2}$-error. Moreover the dynamics of these parameters is itself well-captured by a constant-coefficient system, an averaged modulation system. To some extent, the latter system may be derived formally through an adapted form of geometrical optics. For more heuristics on the latter we refer to \cite{R,R_Roscoff} whereas the corresponding detailed analysis may be found in \cite{JFA-R}. 

Actually, in order to prove Theorem~\ref{thm2} a weaker form of dispersive spectral stability would be sufficient and we have chosen to stick to the stronger form only because it is useful to derive the modulation behavior discussed above and does hold in the cases that we have analyzed. We refer to \cite{JFA-R} for more optimal assumptions. 

\medskip

In the next section, we provide some background and elementary preliminaries to the proof of our main results. In Section~\ref{s:spec-stab} we complete the proof of Theorem~\ref{thm1} by carrying out a small-amplitude/slow/side-band analysis. In Section~\ref{s:linear-stab}  we complete the proof of Theorem~\ref{thm2} by providing a detailed high-frequency analysis, replacing the one of \cite{JFA-R}. 

\medskip

\noindent {\bf Acknowledgment:} L.M.R. expresses his gratitude to INSA Toulouse and C.A. and C.S. theirs to Universit\'e de Rennes 1 for their hospitality during part of the preparation of the present contribution. L.M.R. would like to thank the Isaac Newton Institute for Mathematical Sciences, Cambridge, for support and hospitality during the programme \emph{Dispersive Hydrodynamics}.

\section{Preliminaries}\label{s:pre}

\subsection*{Wave parameterization and small-amplitude limit}

Up to a phase shift, periodic wave profiles may be identified with connected components that are reduced to closed loops of level sets defined by 
\begin{align*}
    \f{1}{2}(\p_x\uu)^2+\mathcal{W}(\uu; c, \lambda)=\mu 
\end{align*}
with potential
\[
\cW(\uu; c, \lambda)=F(\uu)-\f{c}{2}\uu^2-\lambda\uu,
\]
where $F$ is the anti-derivative of $f$ such that $F(0)=0$ and $(\lambda,\mu)$ are wave parameters complementing the phase speed $c$. Throughout, without losing generality, we quotient the phase shift freedom by imposing that $\uu$ is even. For a full discussion of the sense in which the identification is carried out and justifications of other claims of the present subsection, we refer the reader to any of the devoted sections in \cite{JNS-BNR,Nonlinearity-BMR1,IUMJ-BMR,Nonlinearity-BMR2,SMF-Audiard-Rodrigues,NRS}.

Constant states are then characterized as critical points of $\cW(\cdot;c_0,\lambda_0)$ for some $(c_0,\lambda_0)$. Those generating small-amplitude periodic waves are necessarily local minimizers of $\cW(\cdot;c_0,\lambda_0)$ and we are only interested in non degenerate ones, that is on $\uu_0$ such that $\p_u\cW(\uu_0;c_0,\lambda_0)=0$ and $\p_u^2\cW(\uu_0;c_0,\lambda_0)>0$. 

When considering small-amplitude waves, it is sufficient to hold $(c,\lambda)$ fix to the value $(c_0,\lambda_0)$ of the limiting endstate $\uu_0$ and send $\mu$ to $\cW(\uu_0;c_0,\lambda_0)$ from above. We shall do so and consistently we drop the index ${}_0$ on $(c,\lambda)$ and introduce $\delta$ through
\[
\mu^\delta:=\cW(\uu_0; c,\lambda)+\f{\delta^2}{2}\p_u^2 \cW(\uu_0; c,\lambda)
\]
so as to quantify the limit as $\delta\to0$. Near $\uu_0$ the wave profiles $\uu_\delta$ of parameter $(c,\lambda,\mu^\delta)$ have scaled profiles $\uv_\delta(\cdot):=\uu_\delta(k_\delta^{-1}\,\cdot)$ and (fundamental) wave numbers $k_\delta$ expands as
\begin{align*}
\uv_\delta(x)&\stackrel{\delta\to0^+}{=}\uu_0+\delta \cos(2\pi x)+\delta^2
\f{f''(\uu_0)}{12(f'(\uu_0)-c)}\big(\cos(4\pi x)-3\big)+\cO(\delta^3)\,,\\
k_\delta&\stackrel{\delta\to0^+}{=}\frac{\sqrt{f'(\uu_0)-c}}{2\pi}
+\f{\delta^2}{32\pi\sqrt{f'(\uu_0)-c}}\,\left(
{f'''(\uu_0)}
-\frac{5\,(f''(\uu_0))^2}{3\,\left(f'(\uu_0)-c\right)}\right)
+\cO(\delta^4)\,.
\end{align*}

At some point, we shall need to parameterize small-amplitude waves. We shall use $(\delta,\lambda,c)$ and the implicit phase shift as parameters. For later use, we already note here that 
\[
\p_ck_0\,\p_\lambda\uu_0-
\p_\lambda k_0\,\p_c\uu_0
\,=\,-\frac{1}{4\pi\,(f'(\uu_0)-c)^{\frac32}}\neq0\,.
\]
Note that in the foregoing computation we have taken into account that $\uu_0$ depends on $(\lambda,c)$ through $\p_u\cW(\uu_0;c,\lambda)=0$.

\subsection*{Linearization and Bloch symbols}

When discussing stability it is also convenient to scale period to a fixed value. To do so, we introduce $k$ the wave number of the background wave and $\uv(\cdot):=\uu(k^{-1}\,\cdot)$ its scaled profile. We also replace $u$ solving $\p_t u=\cL u$ where $\cL=-\p_x((\p_x^2+f'(\uu)-c)\,\cdot\,)$ with $v$ such that $v(t,\cdot)=u(t,k^{-1}\,\cdot)$. In terms of $v$, the linearized dynamics obeys
\begin{align*}
\p_t v&=L v\,,&
L&:=-k\,\p_x((k\,^2\p_x^2+f'(\uv)-c)\,\cdot\,)\,.
\end{align*}
The operator $L$ is a differential operator with coefficients of period $1$ and we consider it as acting on $L^2(\mR)$ with domain $H^3(\mR)$.

In order to analyze the structure of the spectrum of $L$, we rely on the Bloch-wave decomposition
\begin{align}\label{def-bloch}
    g(x)=\int_{-\pi}^{\pi} e^{i\xi\,x} \check{g}(\xi,x)\,\d \xi
\end{align}
provided by the Floquet-Bloch transform 
\begin{align}\label{def-Bloch}
  \check{g}(\xi,x)
  &:=\sum_{j\in \mZ} e^{i\,2\pi j\,x } \hat{g}(\xi+2j\pi)
  =\frac{1}{2\pi}\sum_{\ell\in\mZ}\,e^{-i\xi\,(x+\ell)}g(x+\ell)\,\,,
\end{align}
which, up to a $\sqrt{2\pi}$ factor, is an isometry from $L^2(\mR)$ to $L^2((-\pi,\pi);L^2(\mR/\mZ))$. Here and throughout we use the following convention for the Fourier transform
\[
\hat{g}(\xi)
:=\frac{1}{2\pi}\int_{\mR} e^{-i\xi\,x} g(x)\,\d x\,.
\]
Note that $(e^{i\xi\,\cdot} \check{g}(\xi,\cdot))(x+1)=e^{i\xi}\,e^{i\xi\,x} \check{g}(\xi,x)$ and, consistently, $\xi$ is called a Floquet exponent and measure a deviation from $1$-periodicity.

The action of $L$ on a function $g$ defined on $\mR$ may be analyzed through
\[
(Lg)(x)=\int_{-\pi}^{\pi} e^{i\xi\,x} L_\xi(\check{g}(\xi,\cdot))(x)\,\d \xi
\]
where the Bloch symbols $L_\xi:=-k\,(\p_x+i\xi)((k^2\,(\p_x+i\xi)^2+f'(\uv)-c)\,\cdot\,)$ are operators on $L^2(\mR/\mZ)$ with domain $H^3(\mR/\mZ)$. Each $L_\xi$ has compact resolvent thus its spectrum is reduced to a discrete set of eigenvalues with finite multiplicity. A key related fact is the spectral decomposition
\[
\sigma(L)=\bigcup_{\xi\in(-\pi,\pi]}\sigma_{per}(L_{\xi})\,,
\]
where we have marked with the index $_{per}$ the fact that $L_\xi$ acts on functions defined over $\mR/\mZ$.

The Hamiltonian structure implies the following relation between $L_\xi^*$ and its adjoint $L_\xi$
\[
(\p_x+i\xi)\,(L_\xi)^*
\,=\,-L_\xi\,(\p_x+i\xi)\,.
\]
Since $(\p_x+i\xi)$ is boundedly invertible on $L^2(\mR/\mZ)$ when $\xi\neq0$, we deduce that the spectrum is symmetric with respect to the imaginary axis, with respect to both location and multiplicity. A continuity argument extends the latter to $L_0$.

\subsection*{Elementary spectral considerations} To discuss the small-amplitude we denote as $L^\delta_\xi$ the Bloch symbol corresponding to $\uv_\delta$.

Simple Fourier series computations provide a full description of the spectrum at $\delta=0$. The eigenvalues of $L_\xi^0$ are given by $\lambda_j(\xi)$, $j\in\mZ$, where 
\begin{align*}
\lambda_j(\xi)&:=\Lambda(2\pi\,j+\xi)\,,&
\Lambda(\zeta)&:=-k_0\,i\zeta\,(k_0^2\,(i\zeta)^2+f'(\uu_0)-c\,)
=-i\,k_0^3\,\zeta\,(4\pi^2-\zeta^2)\,.
\end{align*}
Moreover eigenvalues of $L_\xi^0$ except for the eigenvalue $0$ of $L_0^0$ that have multiplicity three. Furthermore, the third-order derivative of the eigenvalues of $L_\xi^0$ with respect to $\xi$ never vanishes and the second-order derivative vanishes only for the eigenvalue $0$ of $L_\xi^0$.

Since simple eigenvalues perturb smoothly, this already provides a large part of Theorem~\ref{thm1}. There remains only to discuss three things.
\begin{enumerate}
\item What happens to the spectrum near the origin ? 
\item Why does the spectrum far from the origin remains on the imaginary axis ? 
\item Can the perturbation argument be made uniform with respect to the eigenvalue ?
\end{enumerate}
The first point is the object of the next section. 

Concerning the second point, we observe that a direct perturbative argument shows that provided one excludes a neighborhood of the origin in the spectral plane and\footnote{The uniformity of the smallness with respect to the eigenvalue is discussed in the third point.} takes $\delta$ sufficiently small, each eigenvalue of $L_\xi^\delta$ is contained in a ball symmetric with respect to the imaginary axis, not intersecting the real axis and where there is no other eigenvalue. Then the Hamiltonian symmetry of the spectrum recalled above implies that these eigenvalues are purely imaginary.  

As for the third point, it is sufficient to focus on the high frequency regime and to prove that there exists $\delta_0>0$ such that for any $0\leq \delta\leq \delta_0$ and any $(j,\xi)$, $|j|\geq 2$, for some sufficiently small $r_\delta(j,\xi)>0$ one may define the spectral projector
\[
\f{1}{2i \pi} \oint_{\p (B(\lambda_j(\xi),r_\delta(j,\xi)))} (\lambda I-L_{\xi}^{\delta})^{-1} \d \lambda\,.
\]
Incidentally we recall that the rank of a spectral projector is equal to the sum of the algebraic multiplicities of the associated eigenvalues and that the rank map is continuous on projectors. Since the distance between consecutive eigenvalues scales like $j^2$, by enforcing that $j^{-2}\,r_\delta(j,\xi)$ is sufficiently small (uniformly with respect to everything) one deduces from direct and Fourier series computations that when $|\lambda-\lambda_j(\xi)|=r_\delta(j,\xi)$,
\begin{align*}
\|(\lambda I-L_{\xi}^0)^{-1}\|_{L^2(\mR/\mZ)\to H^1_{per}(\mR/\mZ)}&\lesssim \frac{|j|}{r_\delta(j,\xi)}\,,&
\|L_{\xi}^\delta-L_{\xi}^0\|_{H^1_{per}(\mR/\mZ)\to L^2(\mR/\mZ)}&\lesssim \delta\,.&
\end{align*}
This is sufficient to conclude the argument with $r_\delta(j,\xi)$ of size $|j|$.

We have only sketched the present spectral arguments because they are both elementary and quite classical. However we refer to \cite{JDE-Gallay-Haragus} for a similar, fully worked out analysis.

\begin{rmk}\label{rk:hf}
We point out that the arguments sketched here may also be applied to analyze the high-frequency regime of a wave of arbitrary amplitude. Indeed for a fixed wave (possibly of large amplitude), by taking $r(j,\xi)$ such that $j^{-2}\,r(j,\xi)$ is small but $|j|^{-1}\,r(j,\xi)$ is large --- which is possible when $|j|$ is sufficiently large --- and appealing to the Hamiltonian symmetry of the spectrum one deduces that, for any wave, outside a sufficiently large ball of the complex plane the spectrum of $L_\xi$ is made of purely imaginary simple eigenvalues. 
\end{rmk}

\section{Spectral stability -- proof of Theorem \ref{thm1}}\label{s:spec-stab}

To conclude the proof of Theorem~\ref{thm1}, we only need to consider the spectrum of $L_\xi^\delta$ near the origin when $(\delta,\xi)$ is small. Our strategy follows very closely similar computations carried out in \cite{NRS}.

To begin with, we introduce the spectral projector 
\begin{align}\label{projection-0}
    \Pi_{\xi}^{\delta}:=\f{1}{2i \pi} \oint_{S}
    (\lambda I-L_{\xi}^{\delta})^{-1} \d \lambda
\end{align}
where $S$ is a circle centered at $0$ with small radius (independent of $(\delta,\xi)$), which is well-defined when $(\delta,\xi)$ is small. Once a smooth basis $(q_1^{\delta}(\xi, \cdot), q_2^{\delta}(\xi, \cdot), q_3^{\delta}(\xi, \cdot))$ of $\Ran(\Pi_\xi^\delta)$ and a smooth dual basis $(\tq_1^{\delta}(\xi, \cdot), \tq_2^{\delta}(\xi, \cdot), \tq_3^{\delta}(\xi, \cdot))$ of $\Ran((\Pi_\xi^\delta)^*)$ have been chosen, the question is reduced to the study of the spectrum of 
\[
D_{\xi}^{\delta}:=\left(\inp{\tq_j^{\delta}(\xi,\cdot)}{L_{\xi}^{\delta}q_\ell^{\delta}(\xi,\cdot)}_{L_{per}^2}\right)_{1\leq j,\,\ell\leq 3}
\] 
where $\inp{\,\cdot\,}{\,\cdot\,}_{L_{per}^2}$ denotes the canonical scalar product of $L^2(\mR/\mZ)$ skew-linear in its first argument.

\subsection*{Simple computations}

We begin with co-periodic computations, that is, by computing with $\xi=0$. By differentiating the profile equation and using small-amplitude expansions to check independence, one obtains that 
a smooth licit choice when $\xi=0$ is
\begin{align}
q_1^{\delta}(0,\cdot)&:=-\f{\p_x\uvd}{2\pi \delta}\,,\nonumber\\\label{defbases-0}
q_2^{\delta}(0,\cdot)&:=\p_{\delta}\uvd-\p_{\delta} k_{\delta}\f{\p_\lambda\uu_0\p_{c}\uvd
-\p_c\uu_0\p_{\lambda}\uvd}{\p_\lambda\uu_0\p_{c}k_{\delta}
-\p_c\uu_0\p_{\lambda}k_{\delta}}\,,\\
q_3^{\delta}(0,\cdot)&:=\frac{
\p_{c}k_{\delta}\,\p_{\lambda} \uvd -\p_{\lambda} k_{\delta}\p_{c}\uvd}{\p_\lambda\uu_0\p_{c}k_{\delta}
-\p_c\uu_0\p_{\lambda}k_{\delta}}\,,\nonumber
\end{align}
which corresponds to 
\begin{align*}
L_0^{\delta}q_1^{\delta}(0, \cdot)&\equiv 0\,,&
L_0^{\delta}q_2^{\delta}(0, \cdot)&=
E_{12}^{\delta}\delta^2\ 
q_1^{\delta}(0, \cdot)\,,&
L_0^{\delta}q_3^{\delta}(0, \cdot)&=
E_{13}^{\delta}\delta\ q_1^{\delta}(0,\cdot)\,,
\end{align*}
where
\begin{align*}
E_{12}^{\delta}&:= 
-\f{2\pi k_{\delta}\,\p_\lambda\uu_0}{\p_\lambda\uu_0\p_{c}k_{\delta}-\p_c\uu_0\p_{\lambda}k_{\delta}}\f{\p_{\delta} k_{\delta}}{\delta}\,,&
E_{13}^{\delta}&:=-\f{2\pi k_{\delta}\p_{\lambda} k_{\delta}}{\p_\lambda\uu_0\p_{c}k_{\delta}
-\p_c\uu_0\p_{\lambda}k_{\delta}}\,.
\end{align*}
Without even making explicit the dual basis, this choice yields 
\begin{align}\label{D0delta}
    D_{0}^{\delta}&=\left(\begin{array}{ccc}
        0 & E_{12}^{\delta}\delta^2  &  E_{13}^{\delta}\delta \\
       0  & 0 & 0 \\
       0 & 0 &  0 
    \end{array}\right)
\,.
\end{align}

From this choice, we build the basis of the general case through Kato's extension process detailed in \cite[p.199]{Kato}. Explicitly, we set $q_j^{\delta}(\xi, \cdot)=\cU_{\xi}^{\delta}q_j^{\delta}(0, \cdot)$, $j=1,2,3$, where $\cU_{\xi}^\delta$ is defined as the unique solution to the linear Cauchy problem 
\begin{align}\label{def-extop}
\p_\xi\cU_{\xi}^{\delta} 
&=[\p_\xi\Pi_{\xi}^{\delta}, \Pi_{\xi}^{\delta}]\,\cU_{\xi}^{\delta}\,,&
\cU_0^{\delta}&=I.
\end{align}
Similar constructions may be carried out for dual bases but we don't need to make it explicit here.

As scalar constant-coefficient operators, any pair of $L_\xi^0$ and $L_\zeta^0$ commute. This implies that $\p_\xi\Pi_{\xi}^0$ and $\Pi_{\xi}^0$ also commute and thus that $\cU_\xi^0\equiv I$ and $q_j^{0}(\xi, \cdot)=q_j^{0}(0, \cdot)$, $j=1,2,3$. More explicitly 
\begin{align*}
q_1^0 (\xi,\cdot)&=\sin(2\pi \cdot)\,,&
q_2^0(\xi, \cdot)&=\cos (2\pi \cdot)\,,&
q_3^0(\xi, \cdot)&\equiv
\p_{c}k_{\delta}\p_{\lambda} \uu_0 -\p_{\lambda} k_{\delta}\p_{c}\uu_0\,.
\end{align*}
Combining this with the already known trigonometric diagonalization of $L_\xi^0$ we deduce
\begin{align}\label{Dxi0}
\hspace{-1em}
D_{\xi}^{0}&=\begin{pmatrix}
\frac{\lambda_1(\xi)+\lambda_{-1}(\xi)}{2}
&i\frac{\lambda_1(\xi)-\lambda_{-1}(\xi)}{2}&0\\[0.25em]
\frac{\lambda_1(\xi)-\lambda_{-1}(\xi)}{2i}
&\frac{\lambda_1(\xi)+\lambda_{-1}(\xi)}{2}&0\\
0&0&\hspace{-1em}\lambda_0(\xi) 
\end{pmatrix}
=\begin{pmatrix}
ik_0^3\xi (8\pi^2+\xi^2)
&-6\pi k_0^3\xi^2 &0\\
6\pi k_0^3\xi^2 
&\hspace{-1em}
ik_0^3\xi (8\pi^2+\xi^2)&0\\
0&0&\hspace{-1em}
-ik_0^3\xi(4\pi^2-\xi^2)
\end{pmatrix}.
\end{align}

\subsection*{Symmetries and structure}

We now turn to pieces of information that require no computation but rely heavily on the structure. 

As a preliminary we define a bounded operator $\cA$ by $(\cA g)(x)=\overline{g(-x)}$. From the evenness of $\uvd$ stems $L_{\xi}^{\delta} \cA=-\cA \lxd$ and thus 
$\pxd \cA=-\cA\pxd$, and $\cU_{\xi}^{\delta} \cA=\cA\,\cU_{\xi}^{\delta}$. From the latter and the definition \eqref{defbases-0} follows 
\begin{align*}
\cA q_1^{\delta}(\xi,\cdot)&=-q_1^{\delta}(\xi,\cdot)\,,&
\cA q_\ell^{\delta}(\xi,\cdot)&= q_\ell^{\delta}(\xi,\cdot)\,,\quad \ell=2,3\,.
\end{align*}
Jointly with $L_{\xi}^{\delta} \cA=-\cA \lxd$, this implies
\beqs 
\diag (-1,1,1)\,\overline{D_{\xi}^{\delta}}=-D_{\xi}^{\delta}\,\diag (-1,1,1).
\eeqs
thus 
\begin{align*} 
(D_{\xi}^{\delta})_{1\ell}&\in \mR\,,\quad \ell=2,3,&
(D_{\xi}^{\delta})_{j1}&\in \mR\,,\quad j=2,3,&\\
(D_{\xi}^{\delta})_{11}&\in i\mR\,,&
(D_{\xi}^{\delta})_{j\ell}&\in i\mR\,,\quad (j,\ell)\in \{2,3\}^2\,.
\end{align*}

Now we use these pieces of information to gain some knowledge about terms mixing $\delta$ and $\xi$. From the $\xi=0$ computations we already know that $L_0^\delta q_1^{\delta}(0,\cdot)\equiv0$ and, for $j=2,3$, $(L_0^\delta)^*q_j^{\delta}(0,\cdot)\equiv0$. This implies that, for $j=2,3$, 
\[
(\p_\xi D_0^\delta)_{j1}
\,=\,\inp{\tq_j^{\delta}(0,\cdot)}{(\p_\xi L_{0}^{\delta})q_1^{\delta}(0,\cdot)}_{L_{per}^2}\,.
\]
Since the right-hand side of the latter inequality is purely imaginary it must be zero. Likewise since $L_0^\delta q_2^{\delta}(0,\cdot)=\mathcal{O}(\delta^2)$, $(L_0^\delta)^*q_1^{\delta}(0,\cdot)=\mathcal{O}(\delta)$ and $\p_\xi q_2^{\delta}(0,\cdot)=\mathcal{O}(\delta)$, one obtains that
\[
(\p_\xi D_0^0)_{12}
\,=\,\inp{\tq_1^{0}(0,\cdot)}{(\p_\xi L_{0}^{0})q_2^{0}(0,\cdot)}_{L_{per}^2}
\]
is both real and purely imaginary, thus zero.

Gathering various pieces of information one deduces that
\beq\label{def-dxd}
D_\xi^\delta
=\begin{pmatrix}
i\xi B_{1,1}^{\xi,\delta}&B_{1,2}^{\xi,\delta}&\delta\,B_{1,3}^{\xi,\delta}\\[0.5em]
(i\xi)^2\,B_{2,1}^{\xi,\delta}&i\xi\,B_{2,2}^{\xi,\delta}&i\xi\,\delta\,B_{2,3}^{\xi,\delta}\\[0.5em]
(i\xi)^2\delta B_{3,1}^{\xi,\delta}&i\xi\,\delta B_{3,2}^{\xi,\delta}&i\xi\,B_{3,3}^{\xi,\delta}
\end{pmatrix}
\eeq
for some smooth real-valued functions $B_{j,\ell}^{\xi,\delta}$, $1\leq j,\,\ell\leq 3$, with 
\[
B_{1,2}^{\xi,\delta}
\,=\,\delta^2\,\frac12\,\p_\delta^2B_{1,2}^{0,0}
+\,\xi^2\,\frac12\,\p_\xi^2B_{1,2}^{0,0}
+\cO\left(\delta^3+|\xi|^3\right)\,.
\]
This suggests the following blow-up
\[
\widetilde{D}_{\xi}^{\delta}
:=\f{1}{i\xi}P_{\xi}D_{\xi}^{\delta}P_{\xi}^{-1}
=\begin{pmatrix}
B_{1,1}^{\xi,\delta}&B_{1,2}^{\xi,\delta}&\delta\,B_{1,3}^{\xi,\delta}\\[0.5em]
B_{2,1}^{\xi,\delta}&B_{2,2}^{\xi,\delta}&\delta\,B_{2,3}^{\xi,\delta}\\[0.5em]
\delta B_{3,1}^{\xi,\delta}&\delta B_{3,2}^{\xi,\delta}&\,B_{3,3}^{\xi,\delta}
\end{pmatrix}
\]
where $P_{\xi}:=\diag(i\xi,1,1)$. The matrix $\widetilde{D}_{\xi}^{\delta}$ is still smooth and its eigenvalues differ from the eigenvalues of interest by a factor $i\xi$. In particular to control the third-order derivatives, with respect to $\xi$, of the latter we shall study the second-order derivatives of the former. Note moreover that 
\[
\widetilde{D}_{0}^{0}
\,=\,
\begin{pmatrix}
k_0^3\,8\pi^2&0&0\\
-6\pi\,k_0^3&k_0^3\,8\pi^2&0\\
0&0&-k_0^3\,4\pi^2
\end{pmatrix}\,.
\]

\subsection*{Proof of Theorem~\ref{thm1}}

Motivated by the structure of $\widetilde{D}^0_0$, we introduce real smooth $\alpha_1^{\xi,\delta}$, $\alpha_2^{\xi,\delta}$ such that
\begin{align*}
Q_{\xi,\delta} \widetilde{D}_{\xi}^{\delta}Q_{\xi,\delta}^{-1}
&=\begin{pmatrix}
B_{1,1}^{\xi,\delta}+\delta^2\,\alpha_1^{\xi,\delta}B_{1,3}^{\xi,\delta}
&B_{1,2}^{\xi,\delta}+\delta^2\,\alpha_2^{\xi,\delta}B_{1,3}^{\xi,\delta}
&\delta\,B_{1,3}^{\xi,\delta}\\[0.5em]
B_{2,1}^{\xi,\delta}+\delta^2\,\alpha_1^{\xi,\delta}B_{2,3}^{\xi,\delta}
&B_{2,2}^{\xi,\delta}+\delta^2\,\alpha_2^{\xi,\delta}B_{2,3}^{\xi,\delta}
&\delta\,B_{2,3}^{\xi,\delta}\\[0.5em]
0&0
&\,B_{3,3}^{\xi,\delta}-\delta^2\,\left(\alpha_1^{\xi,\delta}B_{1,3}^{\xi,\delta}
+\alpha_2^{\xi,\delta}B_{2,3}^{\xi,\delta}\right)
\end{pmatrix}\,,
\end{align*}
where
\begin{align*}
Q_{\xi,\delta}\,:=\,\begin{pmatrix}
1&0&0\\
0&1&0\\
-\delta\,\alpha_1^{\xi,\delta}&-\delta\,\alpha_2^{\xi,\delta}&0
\end{pmatrix}\,.
\end{align*}
The existence of $\alpha_1^{\xi,\delta}$, $\alpha_2^{\xi,\delta}$ follows from the Inverse Mapping Theorem applied to
\begin{align*}
\begin{pmatrix}\alpha_{1}^{\xi,\delta}&\alpha_{2}^{\xi,\delta}\end{pmatrix}
&\left(\begin{pmatrix}
B_{1,1}^{\xi,\delta}&B_{1,2}^{\xi,\delta}\\[0.5em]
B_{2,1}^{\xi,\delta}&B_{2,2}^{\xi,\delta}\end{pmatrix}
-\,B_{3,3}^{\xi,\delta}\,I_2
\right)\\
&=\,\begin{pmatrix}B_{3,1}^{\xi,\delta}&B_{3,2}^{\xi,\delta}\end{pmatrix}
-\delta^2\,
\begin{pmatrix}
\alpha_1^{\xi,\delta}\left(
\alpha_1^{\xi,\delta}B_{1,3}^{\xi,\delta}
+\alpha_2^{\xi,\delta}B_{2,3}^{\xi,\delta}\right)
&\alpha_2^{\xi,\delta}
\left(\alpha_1^{\xi,\delta}B_{1,3}^{\xi,\delta}
+\alpha_2^{\xi,\delta}B_{2,3}^{\xi,\delta}\right)
\end{pmatrix}\,.
\end{align*}

One of the eigenvalues of $\widetilde{D}_\xi^\delta$ is merely given by
\[
B_{3,3}^{\xi,\delta}-\delta^2\,\left(\alpha_1^{\xi,\delta}B_{1,3}^{\xi,\delta}
+\alpha_2^{\xi,\delta}B_{2,3}^{\xi,\delta}\right)\,.
\]
It is always real and its second-order derivative with respect to $\xi$ is indeed nonzero (since its limiting value is $2\,k_0^3$), as expected.

The two other eigenvalues of $\widetilde{D}_\xi^\delta$ are given as
\[
\frac12\,\left(B_{1,1}^{\xi,\delta}+B_{2,2}^{\xi,\delta}+\delta^2\,\left(\alpha_1^{\xi,\delta}B_{1,3}^{\xi,\delta}+\alpha_2^{\xi,\delta}B_{2,3}^{\xi,\delta}\right)\right)
\pm \sqrt{\Delta^{\xi,\delta}}
\]
where
\[
\Delta^{\xi,\delta}:=\frac14\,\left(B_{1,1}^{\xi,\delta}-B_{2,2}^{\xi,\delta}+\delta^2\,\left(\alpha_1^{\xi,\delta}B_{1,3}^{\xi,\delta}-\alpha_2^{\xi,\delta}B_{2,3}^{\xi,\delta}\right)\right)^2
+\left(B_{1,2}^{\xi,\delta}+\delta^2\,\alpha_2^{\xi,\delta}B_{1,3}^{\xi,\delta}
\right)\left(
B_{2,1}^{\xi,\delta}+\delta^2\,\alpha_1^{\xi,\delta}B_{2,3}^{\xi,\delta}\right)\,.
\]
Since 
\[
\Delta^{\xi,\delta}
\,=\,\delta^2\,\left(\frac14\,\left(\p_\delta\left(B_{1,1}-B_{2,2}\right)^{0,0}\right)^2
+\frac12\,B_{2,1}^{0,0}\,\left(\p_\delta^2B_{1,2}^{0,0}
+2\alpha_2^{0,0}B_{1,3}^{0,0}\right)\right)
+\frac12\,\xi^2\,B_{2,1}^{0,0}\,\p_\xi^2B_{1,2}^{0,0}
+\cO\left(\delta^3+|\xi|^3\right)
\]
with $B_{2,1}^{0,0}\,\p_\xi^2B_{1,2}^{0,0}>0$, one concludes that when $(\p_\delta\left(B_{1,1}-B_{2,2}\right)^{0,0})^2
+2\,B_{2,1}^{0,0}\,(\p_\delta^2B_{1,2}^{0,0}+2\alpha_2^{0,0}B_{1,3}^{0,0})$ is positive small amplitude waves are spectrally stable whereas they are spectrally unstable when it is negative. At this stage we may simply borrow from the related computations in\footnote{We warn the reader that we correct here a sign mistake introduced there at the very end when rewriting in terms of the flux function a correct formula originally expressed in terms of profile potentials.} \cite[Appendix~A.1]{Nonlinearity-BMR2} that $(\p_\delta\left(B_{1,1}-B_{2,2}\right)^{0,0})^2
+2\,B_{2,1}^{0,0}\,(\p_\delta^2B_{1,2}^{0,0}+2\alpha_2^{0,0}B_{1,3}^{0,0})$ is a positive\footnote{Actually it is equal to $k_0^2\Delta_{BF}$ but this does not follow readily from the computations of \cite{Nonlinearity-BMR2}.} multiple of $\Delta_{BF}$. Alternatively one may compute the quantity of interest along the lines of the computations carried out above but, since this computation is rather long and not very instructive, we only provide in Appendix~\ref{s:compute} some intermediate steps for the convenience of the reader wishing to carry it out.

To conclude we observe when $\Delta_{BF}>0$ that, on one hand, though one cannot get joint regularity in $(\delta,\xi)$, one does get for any fixed $\delta>0$, that the two eigenvalues are simple and thus everything is smooth in $\xi$ and that, on the other hand, the second-order derivatives with respect to $\xi$ of the eigenvalues are equivalent to
\[
\mp\frac18\frac{B_{2,1}^{0,0}\,\p_\xi^2B_{1,2}^{0,0}}{\left(\Delta^{\xi,\delta}\right)^{\frac32}}\,,
\]
as $(\delta,\xi)\to(0,0)$, and thus non zero provided that $(\delta,\xi)$ is sufficiently small.

This achieves the proof of Theorem~\ref{thm1}.

\section{Linear stability -- proof of Theorem~\ref{thm2}}\label{s:linear-stab}

We now prove Theorem~\ref{thm2}. As pointed out in the introduction we only need to modify the large eigenvalue expansions of \cite{JFA-R}, that is, to replace \cite[Section~2.4]{JFA-R}. 

\subsection*{Large eigenvalue expansions} 

As in Remark~\ref{rk:hf}, we first point out that outside a sufficiently large neighborhood of the origin, the spectrum of $L_\xi$ is made of simple purely imaginary eigenvalues and they may be parametrized as $(\lambda_j(\xi))_{j\in\mZ,|j|\geq j_0}$ for some large $j_0\in\mN$ in such a way that 
\[
\lambda_j(\xi)\,\stackrel{|j|\to\infty}{=}\,
i\,k^3\,(2\pi\,j+\xi)^3+\cO(|j|)\,,
\]
and $\lambda_j(\xi)=-\lambda_{-j}(-\xi)$, and, for any $\ell\in\mN$, $\lambda_{j}^{(\ell)}(\pi)=\lambda_{j+1}^{(\ell)}(-\pi)$. Likewise, by using spectral projectors and normalizing eigenfunctions, one obtains corresponding smooth direct and dual eigenfunctions $\phi_j(\xi,\cdot)$ and $\tphi_j(\xi,\cdot)$, such that 
\begin{align*}
\|\phi_j(\xi,\cdot)\|_{L^2(\mR/\mZ)}&=\|\tphi_j(\xi,\cdot)\|_{L^2(\mR/\mZ)}\,,&
\langle \tphi_j(\xi,\cdot),\phi_j(\xi,\cdot)\rangle_{L^2_{per}}&=1\,,
\end{align*}
$\overline{\phi_j(\xi,\cdot)}=\phi_{-j}(-\xi,\cdot)$, $\overline{\tphi_j(\xi,\cdot)}=\tphi_{-j}(-\xi,\cdot)$ and, for any $\ell\in\mN$, $\p_\xi^{(\ell)}\phi_{j}(\pi,\cdot)=\p_\xi^{(\ell)}\phi_{j+1}(-\pi,\cdot)$ and $\p_\xi^{(\ell)}\tphi_{j}(\pi,\cdot)=\p_\xi^{(\ell)}\tphi_{j+1}(-\pi,\cdot)$. 

The following proposition gathers the pieces of information that need to be proved on direct and dual eigenfunctions. For the sake of concision, we denote there by $\mZ_{j_0}$ the sets of $j\in\mZ$ such that $|j|\geq j_0$.

\begin{prop}\label{p:hf}  There exist two families of smooth functions over $\mR$ of period one $(r_\ell)_{\ell\in\mN^*}$ and $(\tilde{r}_\ell)_{\ell\in\mN^*}$ and two families of coefficients $(b_\ell(j,\xi))_{(\ell,j,\xi)\in \mN^*\times\mZ_{j_0}\times[-\pi,\pi]}$ and $(\tilde{b}_\ell(j,\xi))_{(\ell,j,\xi)\in \mN^*\times\mZ_{j_0}\times[-\pi,\pi]}$, smooth in $\xi$, such that the following holds. 
\begin{enumerate}
\item For any $\ell\in\mN^*$, any $\alpha\in\mN$, uniformly in $\xi$
\begin{align*}
\p_\xi^\alpha b_\ell(j,\xi)\,&\stackrel{|j|\to\infty}{=}\cO\left(|j|^{-(\ell+\alpha)}\right)\,,&
\p_\xi^\alpha \tilde{b}_\ell(j,\xi)\,&\stackrel{|j|\to\infty}{=}\cO\left(|j|^{-(\ell+\alpha)}\right)\,.
\end{align*}
\item The direct and dual leading order coincide, that is, $r_1=\tilde{r}_1$ and $\tilde{b}_1=b_1$.
\item For any $M\in\mN^*$ and $N\in\mN^*$, one may choose direct and dual eigenfunctions so as to enforce that the remainders
\begin{align*}
R_M(j,\xi,x)&:=\phi_j(\xi,x)
-e^{2i\pi\,j\,x}\left(1+\sum_{\ell=1}^{M}b_\ell(j,\xi)\,r_\ell(x)\right)\,,\\
\tilde{R}_M(j,\xi,x)&:=\tilde{\phi}_j(\xi,x)
-e^{2i\pi\,j\,x}\left(1+\sum_{\ell=1}^{M}\tilde{b}_\ell(j,\xi)\, \tilde{r}_\ell(x)\right)\,,
\end{align*}
satisfy for any $\alpha\in\mN$, $\beta\in\mN$ such that $\alpha\leq N$, $\beta\leq N$, uniformly in $\xi$
\begin{align*}
\p_\xi^\alpha\p_x^\beta R_M(j,\xi,x)\,&\stackrel{|j|\to\infty}{=}\cO\left(|j|^{-(M+1+\alpha-\beta)}\right),\,\quad 
\p_\xi^\alpha\p_x^\beta \tilde{R}_M(j,\xi,x)\,\stackrel{|j|\to\infty}{=}\cO\left(|j|^{-(M+1+\alpha-\beta)}\right)\,.
\end{align*}
\end{enumerate}
\end{prop}

In \cite[Section~2.4]{JFA-R}, there are various lemmas and propositions but their conclusions follow readily from Proposition~\ref{p:hf}

\subsection*{Proof of Proposition~\ref{p:hf}}

To prove Proposition~\ref{p:hf}, it is sufficient to first build direct and dual eigenfunctions $p_j(\xi,\cdot)$ and $\tilde{p}_j(\xi,\cdot)$ with the suitable form of expansion, starting as $e^{2i\pi\,j\,\cdot}$, and then to normalize those into $\phi_j(\xi,\cdot)$ and $\tilde{\phi}_j(\xi,\cdot)$ so as to enforce duality. Since the constructions are similar, we focus only on the analysis of direct eigenfunctions. Note that the fact that the two constructions coincide at order $\cO(|j|^{-2})$ is due to the fact that, when $\lambda\in i\mR$, direct and dual eigenfunction equations themselves coincide up to terms of order $0$, reducing at this order to kernel equations for $\lambda+k^3\,(\p_x+i\xi)^3+(f'(\uv)-c)(\p_x+i\xi)$.

Instead of relying on a direct Bloch-wave analysis as we have done so far, at this stage we turn to a spatial dynamics interpretation of the eigenfunction equation. The equation for a $p_j(\xi,\cdot)$ of period $1$ is equivalent to the fact that $V$ defined by
\[
V_\lambda(x):=e^{i\xi\,x}\begin{pmatrix} p_j\\(\p_x+i\xi)\,p_j\\(\p_x+i\xi)^2p_j
\end{pmatrix}(\xi,x)
\]
satisfies
\begin{align*}
\frac{\d V_\lambda}{\d x}&\,=\,A_\lambda\,V_\lambda\,,&
V_\lambda(1)&=e^{i\,\xi}\,V_\lambda(0)\,,
\end{align*}
where
\[
A_\lambda(x):=\begin{pmatrix}
    0 & 1 & 0  \\
    0 & 0 & 1 \\
    -\f{1}{k^3}\lambda-\f{1}{k^2}\p_x(f'(\uv)) & -\f{1}{k^2}(f'(\uv)-c) & 0
\end{pmatrix}\,.
\]
Note that in this point of view, instead of fixing $\xi$ and looking for $(\lambda,p)$, one fixes $\lambda$ --- here such that $\lambda\in i\mR$ and $|\lambda|$ is sufficiently large --- and looks for $(\xi,V)$.

We shall conclude the proof of Proposition~\ref{p:hf} by designing a suitable fixed point argument for the construction of $(\xi,V)$. To achieve this, we want to bring the ODE system associated with $A_\lambda$ into a diagonal constant-coefficient form up to arbitrarily small remainders. To begin with, we set
\begin{align*}
\omega&:=e^{\frac{2i\pi}{3}}\,,&
\eta&:=
-ik(-i\lambda)^{-\f{1}{3}}\,,&
P&=\diag(1,\eta^{-1},\eta^{-2})\begin{pmatrix}
1&1&1\\
1&\omega&\omega^2\\
1&\omega^2&\omega
\end{pmatrix}\,,
\end{align*}
and observe that for $W_0:=P^{-1}\,V$, the ODE spectral system becomes
\begin{align*}
\frac{\d W_0}{\d x}&\,=\,\left(\eta^{-1}\,B_0+\cR_0\right)\,W_0\,,&
\end{align*}
where $\cR_0$ is a smooth function of $(\eta,x)$, periodic of period $1$ with respect to $x$ and 
\begin{align*}
B_0&:=\diag(1,\omega,\omega^2)\,,&
\cR_0&\stackrel{|\eta|\to0}{=}\cO(1)\,.
\end{align*}
Note that $\eta\in i\mR$ since $\lambda$ is purely imaginary. All the functions we build are smooth in $\eta$ and in the end expansions in Proposition~\ref{p:hf} may be thought as Taylor expansions with respect to this variable. For the sake of readability, most of the times we leave the dependency on $\eta$ implicit.

From here, the transformation to an approximate constant-coefficient diagonal form follows by joining together two classical arguments. The approximate diagonalization is due to the separation of eigenvalues of the leading-order $B_0$ and is comparable to many high-frequency analyses ; we refer for instance to those in \cite{BJNRZ,JNRZ,BJNRZ,JNRYZ,BR,FaR}. In turn the conjugation to constant-coefficient problems stems from the original Floquet theory for periodic ODEs and we need it here in its most trivial form, the one for diagonal systems, or, equivalently, for scalar equations.

Indeed, one may build recursively $Q_M$, $B_M$, $\cR_M$, $M\in\mN$, smooth in $(\eta,x)$ and of period $1$ in $x$, such that for $W_M:=(I+\eta^{M}Q_{M-1})\,(I+\eta^{M-1}Q_{M-2})\cdots(I+\eta Q_0)W_0$ the ODE spectral system becomes
\begin{align}\label{eq:normal-form}
\frac{\d W_M}{\d x}&\,=\,\left(\eta^{-1}\,B_M+\eta^{M}\cR_M\right)\,W_M\,,&
\end{align}
where $B_M$ is constant in $x$ and diagonal. The off-diagonal part of $Q_M$ is determined by requiring that $Q_M\,B_M-B_M\,Q_M+\cR_M$ is diagonal whereas its diagonal part is determined by
\begin{align*}
\frac{\d}{\d x}(Q_M)_{\ell,\ell}\,=\,-\left((\cR_M)_{\ell,\ell}-\int_0^1(\cR_M)_{\ell,\ell}\right)\,,
\qquad \ell=1,2,3,
\end{align*}
and the fact that it is mean-free\footnote{This specific normalization ensures uniqueness but plays no particular role.}. The constant diagonal matrix is defined recursively through
\[
B_{M+1}:=B_M+\eta^{M+1}\diag\left(\int_0^1(\cR_M)_{1,1},\int_0^1(\cR_M)_{2,2},\int_0^1(\cR_M)_{3,3}\right)\,.
\]
Note that at any stage $B_M=B_0+\cO(\eta)$ as $|\eta|\to0$ and the diagonal entries of $B_M$ are two-by-two distinct, with uniform separation.

We would like to use this preparation to build a convenient basis of solutions --- thus, equivalently, a solution operator --- to \eqref{eq:normal-form} hence to the original spectral ODE system. As a preliminary we observe that $\eta$ is a positive multiple of $-i$ and $\Re(i\,\omega)<0<\Re(i\,\omega^2)$. By applying a fixed-point argument to 
\begin{align*}
\Sigma_{M}(x)&=\diag\left(e^{\frac{x}{\eta}(B_M)_{1,1}},e^{\frac{x}{\eta}(B_M)_{2,2}},e^{-\frac{(1-x)}{\eta}(B_M)_{3,3}}\right)\\
&\quad+\eta^M\int_0^x\diag\left(e^{\frac{(x-y)}{\eta}(B_M)_{1,1}},e^{\frac{(x-y)}{\eta}(B_M)_{2,2}},0\right)\,\cR_M(y)\,\Sigma_{M}(y)\,\d y\\
&\quad-\eta^M\int_x^1\diag\left(0,0,e^{-\frac{(y-x)}{\eta}(B_M)_{3,3}}\right)\,\cR_M(y)\,\Sigma_{M}(y)\,\d y
\end{align*}
one obtains a matrix-valued map $\Sigma_M$ defined on $[0,1]$ such that 
\begin{align*}
\frac{\d \Sigma_M}{\d x}&\,=\,\left(\eta^{-1}\,B_M+\eta^{M}\cR_M\right)\,\Sigma_M\,,
\qquad \textrm{on }[0,1]
\end{align*}
and the remainder $\widetilde{\Sigma}_M$ defined by
\[
\widetilde{\Sigma}_M(x):=\Sigma_{M}(x)-\diag\left(e^{\frac{x}{\eta}(B_M)_{1,1}},e^{\frac{x}{\eta}(B_M)_{2,2}},e^{-\frac{(1-x)}{\eta}(B_M)_{3,3}}\right)
\]
satisfies
\begin{align*}
\widetilde{\Sigma}_M(0)&=\begin{pmatrix}
0&0&0\\0&0&0\\ *&*&*
\end{pmatrix}\,,&
\widetilde{\Sigma}_M(1)&=\begin{pmatrix}
*&*&*\\ *&*&*\\0&0&0
\end{pmatrix}\,,\\
\frac{\d^\ell \widetilde{\Sigma}_M}{\d x^\ell}(x)
&\stackrel{|\eta|\to0}{=}\cO\left(|\eta|^{M-2\ell}\right)\,,&x\in[0,1],\quad
0\leq\ell\leq M\,.
\end{align*}
Note that at this stage we have made a particular choice and traded infinite regularity for finite --- but arbitrarily large --- regularity.

Now we may obtain $e^{i\xi}$ under the form $\exp(\eta^{-1}\,(B_M)_{1,1})\,\gamma_M$ by solving
\[
\det\left(\Sigma_{M}(1)-e^{\eta^{-1}\,(B_M)_{1,1}}\,\gamma_M\,\Sigma_{M}(0)\right)\,=\,0
\]
with $\gamma_M\stackrel{|\eta|\to0}{=}1+\cO(|\eta|^{M})$. This also yields
\begin{align*}
\frac{\d^\ell \gamma_M}{\d \eta^\ell}(\eta)
&\stackrel{|\eta|\to0}{=}\cO\left(|\eta|^{M-\ell}\right)\,,&0\leq\ell\leq \frac{M}{2}\,.
\end{align*}
Since $\gamma_M$ is valued in a small neighborhood of $1$, one may use a standard determination of the logarithm to define $(j,\xi)$ with arbitrarily large regularity through
\[
i(2\pi\,j+\xi):=\eta^{-1}\,(B_M)_{1,1}+\ln(\gamma_M)\,.
\]
Note that our current proof does not yield that the foregoing expression is indeed purely imaginary, we deduce it from our a priori acknowledge that $\lambda$ does lie in the spectrum of $L_\xi$ for some $\xi$ and that there is no other possibility. Likewise one may solve 
\[
\Sigma_{M}(1)W_M^0=e^{\eta^{-1}\,(B_M)_{1,1}}\,\gamma_M\,\Sigma_{M}(0)W_M^0
\]
for $W_M^0$ in a $\cO(|\eta|^{M})$-neighborhood of $(1,0,0)$. This also comes with similar bounds for derivatives with respect to $(\eta,x)$ up to some arbitrarily large order. By setting $W_M(x):=\Sigma_M(x)W_M^0$, and undoing the two changes from $(V_\lambda,\lambda)$ to $(W_M,\eta)$ and from $(p_j(\xi,\cdot),j,\xi)$ to $(V_\lambda,\lambda)$, one obtains the sought $p_j(\xi,\cdot)$. Let us stress that to conclude the proof we simply need to take $M$ sufficiently large and to expand with respect to $\eta$ all the smooth functions involved in the construction.

\appendix

\section{The Benjamin-Feir index}\label{s:compute}

We provide here some of the intermediate steps to prove that
\[(\p_\delta\left(B_{1,1}-B_{2,2}\right)^{0,0})^2
+2\,B_{2,1}^{0,0}\,(\p_\delta^2B_{1,2}^{0,0}+2\alpha_2^{0,0}B_{1,3}^{0,0})
\,=\,k_0^2\Delta_{BF}\,.
\]
From \eqref{D0delta} and \eqref{Dxi0}, we already deduce that
\begin{align*}
B_{21}^{00}&=-6\pi k_0^3\,,&
\p_{\delta}^2 (B_{12})^{00}
&=-\f{1}{24\pi k_0}\left[5 (f''(\uu_0))^2-3 f'''(\uu_0)(2\pi k_0)^2\right]\,,&  
B_{13}^{00}&=-2\pi k_0\, f''(\uu_0)\,.
\end{align*}

Then computing profile expansions, one may obtain in a lengthy but straightforward way that
\begin{align*}
\p_\delta B_{1,1}^{0,0}&=0\,,&\p_\delta B_{2,2}^{0,0}&=0\,.
\end{align*}

At last, the equation defining $(\alpha_1^{0,0},\alpha_2^{0,0})$ yields
\[
\alpha_2^{0,0}\,=\,\frac{B_{3,2}^{0,0}}{B_{2,2}^{0,0}-B_{3,3}^{0,0}}
\,=\,\frac{B_{3,2}^{0,0}}{12\pi\,k_0^3}\,.
\]
In turn, another lengthy but straightforward computation gives 
\[
B_{3,2}^{0,0}\,=\,k_0\,f''(\uu_0)\,.
\]
With this in hand we have all the pieces to conclude the computation.

\bibliographystyle{alphaabbr}
\bibliography{refgkdv}
\end{document}